\documentclass[11pt]{article}
\usepackage{enumerate}
\usepackage{amssymb,a4wide,latexsym,makeidx,epsfig,fleqn}
\usepackage{amsthm}
\usepackage{amsmath}
\usepackage{enumerate}
\newtheorem{theorem}{Theorem}[section]
\newtheorem{conjecture}[theorem]{Conjecture}

\newtheorem{definition}[theorem]{Definition}
\newtheorem{lemma}[theorem]{Lemma}
\newtheorem{proposition}[theorem]{Proposition}

\begin{document}
\textwidth 150mm \textheight 225mm
\title{The largest signless Laplacian spectral radius of uniform supertrees with diameter and pendent edges (vertices)
\thanks{ Supported by the National Natural Science Foundation of China (No. 11871398), the Natural Science Basic Research Plan in Shaanxi Province of China (Program No. 2018JM1032) and the Seed Foundation of Innovation and Creation for Graduate Students in Northwestern Polytechnical University (No. ZZ2018171).}}
\author{{Cunxiang Duan, Ligong Wang\footnote {Corresponding author.} and Peng Xiao}\\
{\small Department of Applied Mathematics, School of Science, Northwestern
Polytechnical University,}\\ {\small Xi'an, Shaanxi 710072,
People's Republic
of China.}\\ {\small E-mail:cxduanmath@163.com;lgwangmath@163.com;xiaopeng@sust.edu.cn}\\
}
\date{}
\maketitle
\begin{center}
\begin{minipage}{120mm}
\vskip 0.3cm
\begin{center}
{\small {\bf Abstract}}
\end{center}
{\small Let $S_{1}(m, d, k)$ be the $k$-uniform supertree obtained from a loose path $P:v_{1}, e_{1}, v_{2}, \ldots,v_{d}, e_{d}, v_{d+1}$ with length $d$ by attaching $m-d$ edges at vertex $v_{\lfloor\frac{d}{2}\rfloor+1}.$ Let $\mathbb{S}(m,d,k)$ be the set of $k$-uniform supertrees with $m$ edges and diameter $d$ and $q(G)$ be the signless Laplacian spectral radius of a $k$-uniform hypergraph $G$. In this paper, we mainly determine $S_{1}(m,d,k)$ with the largest signless Laplacian spectral radius among all supertrees in $\mathbb{S}(m,d,k)$ for $3\leq d\leq m-1$. Furthermore, we determine the unique uniform supertree with the maximum signless Laplacian spectral radius among all the uniform supertrees with $n$ vertices and pendent edges (vertices).

 }
\end{minipage}
\end{center}

\section{Introduction }
\label{sec:ch6-introduction}

A hypergraph $G=(V,E)$ on $n$ vertices is a set of vertices, say $V=V(G)=\{u_{1},u_{2},\ldots,u_{n}\}$ and a set of edges, say $E=E(G)=\{e_{1},e_{2},\ldots,e_{m}\}$. A hypergraph is called $k$-uniform if every edge contains precisely $k$ vertices. A supertree is a hypergraph which is both connected and acyclic \cite{LiSQ}.

For a $k$-uniform hypergraph $G$, a path of length $l$ is defined to be an alternating sequence of vertices and edges $u_{1}, e_{1}, u_{2},\ldots, u_{l},e_{l}, u_{l+1},$ where $u_{1},u_{2},\ldots, u_{l+1}$ are distinct vertices of $G$, $e_{1},e_{2},\ldots,e_{l}$ are distinct edges of $G$ and $u_{i},u_{i+1}\in e_{i}$ for $i=1,2,\ldots,l.$ If there exists a path between any two vertices of $G$, then $G$ is called connected. The distance between two vertices is the length of the shortest path connecting them. The diameter $d$ of a connected $k$-uniform hypergraph $G$ is the maximum distance among all vertices of $G$.

The degree of a vertex $u$ (denoted by $d_{G}(u)$ or $d_u$) of a $k$-uniform hypergraph $G$ is the number of edges containing $u$. For a $k$-uniform hypergraph $G$ with $V(G)=\{u_1, u_2,\ldots,u_n\}$, if $E(G)=\{e_1, e_2, \ldots,e_m\}$ with $m=\frac{n-1}{k-1}$, where $e_i=\{u_{(i-1)(k-1)+1},\ldots, u_{(i-1)(k-1)+k}\}$ for $i=1,2,\ldots,m$, then $k$-uniform hypergraph $G$ is called a $k$-uniform loose path, denoted by $P_{n,k}$. A vertex of degree one in $e_1$ is an end vertex of $P_{n,k}$. A path $P=(u_0,e_1,u_1,\ldots,e_p,u_p)$ in a $k$-uniform hypergraph $G$ is called a pendent path at $u_0$, if $d_{u_0} \geq 2$, $d_{u_i}=2$ for $1 \leq i \leq p-1$, $d_{u}=1$ for $u \in e_{i}\setminus \{u_{i-1},u_{i}\}$ with  $1 \leq i \leq p$, and $d_{u_p}=1$. If $p=1$, then $e_1$ is a pendent edge of $G$. 

A tensor $\mathcal{A}$ with order $k$ and dimension $n$ over the complex field $\mathbb{C}$ is a multidimensional array $$\mathcal{A}=(a_{i_{1}i_{2}\ldots i_{k}}),~1\leq i_{1},i_{2},\ldots,i_{k}\leq n.$$ The tensor $\mathcal{A}$ is called symmetric if its entries are invariant under any permutation of their indices.

The adjacency tensor \cite{CoDu} of a $k$-uniform hypergraph $G$ with $n$ vertices, denoted by $\mathcal{A}(G)$, is an order $k$ dimension $n$ symmetric tensor with entries $$ a_{i_{1}i_{2}\ldots i_{k}}=\left\{
\begin{array}{ll}
\frac{1}{(k-1)!},& \mbox {if}   ~\{i_{1},i_{2},\ldots, i_{k}\} \in E(G),
\\
0,& \mbox {otherwise}.
\end{array}
\right.$$

\noindent\begin{definition}\label{de:ch-1}(\cite{Lim,Q}) Let $\mathcal{A}$ be an order $k$ dimension $n$ tensor, and $x =(x_{1}, x_{2}, \ldots, x_{n})^{T}\in \mathbb{C}^{n}$ be a column vector of dimension $n$. Then $\mathcal{A}x^{k-1}$ is defined to be a vector in $\mathbb{C}^{n}$ whose $i$th component is the following:
$$(\mathcal{A}x^{k-1})_{i}=\sum\limits^{n}_{i_{2},\ldots,i_{k}=1}a_{ii_{2}\ldots i_{k}}x_{i_{2}}x_{i_{3}}\ldots x_{i_{k}}, ~( 1\leq i \leq n).$$
\end{definition}

Let $x^{[k-1]}=(x_{1}^{k-1},x_{2}^{k-1},\ldots,x_{n}^{k-1})^{T}\in \mathbb{C}^{n}.$~If $ \mathcal{A} x^{k-1}=\lambda x^{[k-1]}$ has a solution $x\in \mathbb{C}^{n}\setminus\{\mathbf{0}\}$, then $\lambda$ is called an eigenvalue of $\mathcal{A}$ and $x$ is an eigenvector associated with $\lambda$.

Let $\mathcal{A}$ be an order $k$ dimension $n$ tensor. The spectral radius of $\mathcal{A}$ is defined as
\begin{equation*}
\rho(\mathcal{A})=\max\{\mid \lambda \mid:\lambda \mbox{ is an eigenvalue of } \mathcal{A}\}.
\end{equation*}

We call $\rho(\mathcal{A})$ the spectral radius of a $k$-uniform hypergraph $G$, denoted by $\rho(\mathcal{A})=\rho(G)$.

\noindent\begin{theorem}\label{de:th-2}(\cite{Qi}) Let $\mathcal{A}$ be an $k$ order $n$ dimension nonnegative symmetric tensor. Then we have $$\rho(\mathcal{A})=\max\{x\mathcal{A}x^{k-1}|\sum_{i=1}^{n}x_{i}^{k}=1,x\in\mathbb{R}^{n}_{+}\},$$
where, $\mathbb{R}^{n}_{+}=\{x\in\mathbb{R}^{n}:x\geq0\}.$ Furthermore, $x\in\mathbb{R}^{n}_{+}$ with $\sum_{i=1}^{n}x_{i}^{k}=1$ is an optimal solution of above optimization problem if and only if it is an eigenvector corresponding to the eigenvalue $\rho(\mathcal{A}).$
\end{theorem}

Let $\mathcal{D}=\mathcal{D}(G)$ be a $k$ order $n$ dimension diagonal tensor with its diagonal element $d_{ii\ldots i}$ being $d_{i}$, the degree of vertex $i$, for all $i \in[n]$. Then $\mathcal{Q}(G)=\mathcal{D}(G)+\mathcal{A}(G)$ is the signless Laplacian tensor of the hypergraph $G$ \cite{QiL}. The signless Laplacian eigenvalues refer to the eigenvalues of the signless Laplacian tensor. Let $q(G)$ be the signless Laplacian spectral radius of $G$. If $G$ is a connected $k$-uniform hypergraph, then there exists a unique positive eigenvector $x$ corresponding to $q(G)$ with $\sum_{i=1}^{n}x_{i}^{k}=1.$ Such positive eigenvector is called the principal eigenvector of $\mathcal{Q}(G).$ For a vertex $i\in V,$ we simplify $E_{\{i\}}$ as $E_{i}.$ It is the set of edges containing the vertex $i,$ i.e., $E_{i}=\{e\in E| i\in e\}.$

It is easy to calculate for the signless Laplacian tensor $\mathcal{Q}(G)$ that
\begin{equation*}
 x^{T}(\mathcal{Q}(G)x^{k-1})=\sum\limits_{\{i_{1},i_{2},\ldots,i_{k}\}\in E(G)}(x_{i_{1}}^{k}+x_{i_{2}}^{k}+\cdots+x_{i_{k}}^{k}+kx_{i_{1}}x_{i_{2}}\cdots x_{i_{k}}).
\end{equation*}

The following result can be obtained directly from Definition 1.1 and will be used in the sequel.

\noindent\begin{theorem}\label{de:th-3}(\cite{LiSQ})
Let $\mathcal{Q}(G)$ be an order $k$ dimension $n$ signless Laplacian tensor of the hypergraph $G$ and $x =(x_{1}, x_{2}, \ldots, x_{n})^{T}\in \mathbb{C}^{n}$ be a column vector of dimension $n$. Then $\mathcal{Q}(G)x^{k-1}$ is defined to be a vector in $\mathbb{C}^{n}$ whose $i$th component is the following:
\begin{equation*}
(\mathcal{Q}(G)x^{k-1})_{i}=d_{i}x_{i}^{k-1}+\sum\limits_{e\in E_{i}(G)}x_{i_{2}}x_{i_{3}}\ldots x_{i_{k}}=\sum\limits_{e\in E_{i}(G)}(x_{i}^{k-1}+x_{i_{2}}x_{i_{3}}\ldots x_{i_{k}})~( 1\leq i \leq n).
\end{equation*}
\end{theorem}

Spectral graph theory has a long history behind its development \cite{BrHa,CVRS,GuoS,GuXC}. Guo and Shao \cite{GuoS} determined the first $\lfloor\frac{d}{2}\rfloor+1$ spectral radii trees among all trees with a given fixed diameter. In 2008, Lim \cite{Lim} proposed the study of the spectra of hypergraphs by using the spectra of tensors. Recently, there has been a lot of activity concerning spectral hypergraph theory \cite{LiZB,LiSQ,LimL,LiKY,LuMa,OYQY,YuSS}. Li, Shao and Qi \cite{LiSQ} gave the operations of moving edges, edge-releasing and total grafting on hypergraphs, and determined the first two spectral radii of $k$-uniform supertrees with given $n$ vertices. Many scholars investigated the extreme spectral radius of $k$-uniform hypergraphs on other conditions \cite{LuMa,OYQY,YuSS}. Xiao, Wang and Lu \cite{XWL} investigated the spectral radius of the $k$-uniform hypergraph when the uniform hypergraph is perturbed by 2-switch operation, and determined the maximum spectral radius of $k$-uniform supertrees with given a degree sequence. Kang, Liu and Shan \cite{KLS} presented a lower bound for the spectral radius in terms of vertex degrees and characterized the extremal hypergraphs, and proved a lower bound for the signless Laplacian spectral radius concerning degrees and gave a characterization of the extremal hypergraphs. Xiao, Wang and Du \cite{XiWD} determined the first two largest spectral radii of uniform supertrees with given diameter.  Su et al. \cite{SKLS} determine the first $\lfloor\frac{d}{2}\rfloor+1$ largest spectral radii of $k$-uniform supertrees with size $m$ and diameter $d$ and the first two smallest spectral radii of supertrees with size $m$ by using the methods of grafting operations on supertrees and comparing matching polynomials of supertrees. Xiao and Wang \cite{XW} determined the maximum spectral radius of uniform hypergraphs with given number of pendant edges and pendent vertices, respectively. Many scholars also started to investigate the signless Laplacian spectral radius of $k$-uniform hypergraphs \cite{LMZW,QiL,YZLQ}. 

Let $\mathbb{S}(m,d,k)$ be the set of $k$-uniform supertrees with $m$ edges and diameter $d.$ And let $\mathbb{T}(n,p,k)$ be the set of $k$-uniform supertrees with $n$ vertices and $p$ pendent edges. Let $\mathbb{G}(n,q,k)$ be the set of $k$-uniform supertrees with $n$ vertices and $q$ pendent vertices. In this paper, we determine the supertree with the largest signless Laplacian spectral radius among all supertrees in $\mathbb{S}(m,d,k)$ for $3\leq d\leq m-1$. We also determine the unique supertree with the second largest signless Laplacian spectral radius among all supertrees in $\mathbb{S}(m,3,k)$. And we respectively determine the supertree with the largest signless Laplacian spectral radius among all supertrees in $\mathbb{T}(n,p,k)$ and $\mathbb{G}(n,q,k).$

In Section 2, some necessary notations and lemmas are given. In Section 3, we mainly determine the unique $k$-uniform supertree attains the largest spectral radius in $\mathbb{S}(m,d,k).$ In Section 4, we respectively determine the supertree with the largest signless Laplacian spectral radius among all supertrees in $\mathbb{T}(n,p,k)$ and $\mathbb{G}(n,q,k).$

\section{Preliminaries}
\label{sec:ch-sufficient}

In this section, we give some useful notations and lemmas.

The following definition can be found in \cite{HuQS}.

\noindent\begin{definition}\label{de:ch-1}(\cite{HuQS}) Let $G=(V,E)$ be a graph. For any $k\geq3$, the $k$th power of $G,$ $G^{k}=(V^{k},E^{k})$ is defined as the $k$-uniform hypergraph with edges $E^{k}: = \{e \cup \{i_{e,1},\ldots,i_{e,k-2}\} | e\in E\}$ and the set of vertices $V^{k}: = V \cup ( \bigcup_{e\in E} \{i_{e,1},\ldots,i_{e,k-2}\}).$
\end{definition}

A $k$-uniform hypergraph $G$ is a $k$th power hypergraph (of some graph) if and only if each edge of $G$ contains at least $k-2$ pendent vertices. It is obvious that a $k$-uniform loose path is the $k$th power of a path.

\noindent\begin{definition}\label{de:ch-2}(\cite{HuQS}) The $k$th power of a tree is called a $k$-uniform hypertree.
\end{definition}

The hypertree $S_{n,k}$ is called a hyperstar, which is the $k$th power of an ordinary star $S_{n'}$ with the number of vertices $n=(n'-1)(k-1)+1$ (thus we have $n'=\frac{n-1}{k-1}+1$ here).

Li, Shao and Qi \cite{LiSQ} introduced the operation of moving edges on hypergraphs.

\noindent\begin{definition}\label{de:ch-3}(\cite{LiSQ})
Let $r\geq1,$ $G=(V, E)$ be a $k$-uniform hypergraph with $u \in V$ and $e_{1}, \ldots, e_{r}\in E,$ such that $u \notin e_{i}$ for $i =1, \ldots, r.$ Suppose that $v_{i}\in e_{i}$ (the vertices $v_{1}, \ldots, v_{r}$ need not be distinct) and write $e_{i}'=(e_{i}\setminus\{v_{i}\})\cup \{u\}~(i =1, 2, \ldots, r).$ Let $G'=(V, E')$ be the hypergraph with $E'=(E\setminus\{e_{1}, \ldots, e_{r}\})\cup\{e'_{1}, \ldots, e'_{r}\}.$ Then we say that $G'$ is obtained from $G$ by moving edges $(e_{1},\ldots, e_{r})$ from $(v_{1}, \ldots, v_{r})$ to $u.$
\end{definition}

\noindent\begin{lemma}\label{le:ch-4}(\cite{LiSQ})
Let $r\geq1,$ $G$ be a connected $k$-uniform hypergraph, $G'$ be the hypergraph obtained from $G$ by moving edges $(e_{1}, \ldots, e_{r})$ from $(v_{1}, \ldots, v_{r})$ to $u.$ Let $x$ be the principal eigenvector of $\mathcal{Q}(G)$ corresponding to $q(G).$ Suppose that $x_{u}\geq \max\limits_{1\leq i\leq r}\{x_{v_{i}}\},$ then $q(G') >q(G).$
\end{lemma}

The following edge-releasing operation on hypergraphs is a special case of the above defined moving edge operation.

\noindent\begin{definition}\label{de:ch-5}(\cite{LiSQ})
Let $G$ be a $k$-uniform supertree, $e$ be a non-pendent edge of $G$ and $u\in e.$ Let $e_{1}, e_{2}, \ldots, e_{r}$ be all the edges of $G$ adjacent to $e$ but not containing $u,$ and suppose that $e_{i}\cap e=\{v_{i}\}$ for $i =1,\ldots, r.$ Let $G'$ be the hypergraph obtained from $G$ by moving edges $(e_{1}, e_{2}, \ldots, e_{r})$ from $(v_{1}, v_{2}, \ldots, v_{r})$ to $u.$ Then $G'$ is said to be obtained from $G$ by an edge-releasing operation on $e$ at $u.$
\end{definition}

\noindent\begin{lemma}\label{le:ch-6}(\cite{LiSQ})
Let $G'$ be a supertree obtained from a $k$-uniform supertree $G$ by edge-releasing a non-pendent edge $e$ of $G$ at $v.$ Then $q(G') >q(G).$
\end{lemma}

In the following, Xiao, Wang and Du \cite{XiWD} gave the effection of edge-releasing in the comparison of diameter of $k$-uniform  hypergraphs.

\noindent\begin{lemma}\label{le:ch-7}(\cite{XiWD})
Let $G'$ be a supertree obtained from a k-uniform supertree $G$ by edge-releasing a non-pendent edge $e$ of $G$ at $v.$ Then $d(G')\leq d(G).$
\end{lemma}

In a $k$-uniform hypergraph, an edge $e$ is called a branch edge if $e$ contains at least 3 non-pendent vertices. If $e$ is not a branch edge, then it is called a non-branch edge.

In the following, we will study some operations and its applications in the comparison of the signless Laplacian spectral radius and diameter of hypergraphs.

\noindent\begin{lemma}\label{le:ch-8}
Let $G$ be a $k$-uniform supertree, and $e=\{v_{1}, v_{2},\ldots, v_{k}\}$ be a branch edge of $G$. Suppose that $v_{1}, v_{2}, \ldots, v_{r}$ are all the non-pendent vertices of $e,$ and $e_{i,1}, e_{i,2}, \ldots, e_{i,t_{i}}$ are all the edges (except $e$) incident with $v_{i},$ $3 \leq i\leq r.$ Let $G'$ be the hypergraph obtained from $G$ by moving edges $e_{i,1}, e_{i,2}, \ldots, e_{i,t_{i}}$ from $v_{i}$ to $v_{1},$ $3 \leq i\leq r.$ Then $q(G')>q(G)$ and $d(G')\leq d(G).$
\end{lemma}

\noindent {\bf Proof.} Let $x$ be the principal eigenvector of $\mathcal{Q}(G)$ corresponding to $q(G)$, and $x_{v_{s}}=\max\limits_{3\leq i\leq r}\{x_{v_{i}}\}.$ It is obvious that either $x_{v_{1}}\geq x_{v_{s}}$ or $x_{v_{s}}>x_{v_{1}}$ holds. If $x_{v_{1}}\geq x_{v_{s}},$ by Lemma 2.4, then $q(G') >q(G).$ If $x_{v_{s}}>x_{v_{1}},$ then we obtain a hypergraph $G^{''}$ from $G$ by moving edges $e_{i,1}, e_{i,2}, \ldots, e_{i,t_{i}}$ from $v_{i}$ to $v_{s}$ for $3 \leq i\leq r$, $i \neq s,$ and moving all the edges (except $e$) incident with $v_{1}$ from $v_{1}$ to $v_{s}.$ By Lemma 2.4, we have $q(G^{''}) >q(G).$ It is obvious that $G'\cong G^{''}.$ Hence $q(G') =q(G^{''}) >q(G).$
It is proved that $d(G')\leq d(G)$ in \cite{XiWD}, we omit this proof.
This completes the proof. \hfill$\square$

Xiao, Wang and Lu \cite{XWL} studied the effection of the spectral radius of a hypergraph under 2-switch operation.

Let $G=(V,E)$ be a $k$-uniform hypergraph, and $e=\{u_{1}, u_{2}, \ldots , u_{k}\}, ~f=\{v_{1}, v_{2}, \ldots , v_{k}\}$ be two edges of $G.$ Let $e'=(e\setminus U_{1})\cup V_{1}, ~f'=(f\setminus V_{1})\cup U_{1},$ where $U_{1}=\{u_{1}, u_{2}, \ldots , u_{r}\}, ~V_{1} = \{v_{1}, v_{2}, \ldots , v_{r}\},~ 1 \leq r \leq k-1.$ Let $G'=(V,E')$ be the hypergraph with $E'=(E \setminus\{e, f\})\cup\{e', f'\}.$ Then we say that $G'$ is obtained from $G$ by $e \mathop{\rightleftharpoons}\limits_{v_1,v_2,\ldots,v_r}^{u_1,u_2,\ldots,u_r}f$ or $e \mathop{\rightleftharpoons}\limits_{V_{1}}^{U_{1}} f.$

Let $G = (V,E)$ be a connected $k$-uniform hypergraph, and $x$ be a vector of dimension $n$. For the simplicity of the notation, we write: $x_{V'}=\prod_{v\in V'}x_{v}, ~V'\subseteq V.$

\noindent\begin{lemma}\label{le:ch-9} Let $G =(V,E)$ be a connected $k$-uniform hypergraph, and $e=\{u_{1}, u_{2}, \ldots , u_{k}\},$ $f=\{v_{1}, v_{2}, \ldots , v_{k}\}$ be two edges of $G$ such that $e\cap f=\phi$. Let $G'$ be a connected $k$-uniform hypergraph obtained from $G$ by $e \mathop{\rightleftharpoons}\limits_{V_{1}}^{U_{1}} f,$ where $U_{1} = \{u_{1}, \ldots, u_{r}\},$ $U_{2}=e \setminus U_{1},$ $V_{1} = \{v_{1}, \ldots , v_{r}\},$~$V_{2} = f\setminus V_{1},~1 \leq r \leq k-1.$ Let $x$ be the principal eigenvector of $\mathcal{Q}(G)$. If $x_{U_{1}} \geq x_{V_{1}}$ and $x_{U_{2}}\leq x_{V_{2}},$ then $q(G) \leq q(G').$ Moreover, if one of the two inequations is strict, then $q(G) < q(G').$
\end{lemma}
\noindent {\bf Proof.}
Note that $d_{G}(u)=d_{G'}(u)$ for any $u\in V(G)$.
By Theorem 1.3, we have

\begin{equation*}
\begin{aligned}
q(G')-q(G)&\geq\sum_{u \in V(G)}d_{G'}(u)x_u^{k}+\sum\limits_{\{w_{1},w_{2},\ldots,w_{k}\}=e\in E(G')}(kx_{w_{1}}x_{w_{2}}\ldots x_{w_{k}})\\&- \sum_{u \in V(G)}d_{G}(u)x_u^{k}-\sum\limits_{\{w_{1},w_{2},\ldots,w_{k}\}=e\in E(G)}(kx_{w_{1}}x_{w_{2}}\ldots x_{w_{k}})\\
\end{aligned}
\end{equation*}
\begin{equation*}
\begin{aligned}
&=k(x_{V_{1}}x_{U_{2}}+x_{U_{1}}x_{V_{2}}-x_{U_{1}}x_{U_{2}}-x_{V_{1}}x_{V_{2}})
\\&=k(x_{U_{2}}(x_{V_{1}}-x_{U_{1}})+x_{V_{2}}(x_{U_{1}}-x_{V_{1}}))
\\&=k(x_{V_{2}}-x_{U_{2}})(x_{U_{1}}-x_{V_{1}})\\&\geq0,
\end{aligned}
\end{equation*}

If $q(G')=q(G),$ then we have that $x$ is an eigenvector of $\mathcal{Q}(G')$ corresponding to $q(G').$ Hence $\mathcal{Q}(G)x^{k-1}=q(G)x^{[k-1]}$ and $\mathcal{Q}(G')x^{k-1}=q(G')x^{[k-1]},$ where $\mathcal{Q}(G)$ and $\mathcal{Q}(G')$ are signless Laplacian tensors of $G$ and $G'$, respectively. Then
\begin{equation}
 q(G)x_{v_{1}}^{k-1}=d_{G}(v_1)x_{v_{1}}^{k-1}+x_{V_{1}\backslash\{v_{1}\}}x_{V_{2}}+\sum\limits_{\{v_{1},w_{2},\ldots,w_{k}\}\in E(G)\cap E(G')}(x_{w_{2}}x_{w_{3}} \cdots x_{w_{k}}),
\end{equation}

\begin{equation}
 q(G')x_{v_{1}}^{k-1}=d_{G'}(v_1)x_{v_{1}}^{k-1}+x_{V_{1}\backslash\{v_{1}\}}x_{U_{2}}+\sum\limits_{\{v_{1},w_{2},\ldots,w_{k}\}\in E(G)\cap E(G')}(x_{w_{2}}x_{w_{3}} \cdots x_{w_{k}}),
\end{equation}

\begin{equation}
 q(G)x_{v_{r+1}}^{k-1}=d_{G}(v_1)x_{v_{r+1}}^{k-1}+x_{V_{1}}x_{V_{2}\backslash\{v_{r+1}\}}+\sum\limits_{\{v_{r+1},w_{2},\ldots,w_{k}\}\in E(G)\cap E(G')}(x_{w_{2}}x_{w_{3}} \cdots x_{w_{k}}),
\end{equation}

\begin{equation}
 q(G')x_{v_{r+1}}^{k-1}=d_{G'}(v_1)x_{v_{r+1}}^{k-1}+x_{U_{1}}x_{V_{2}\backslash\{v_{r+1}\}}+\sum\limits_{\{v_{r+1},w_{2},\ldots,w_{k}\}\in E(G)\cap E(G')}(x_{w_{2}}x_{w_{3}} \cdots x_{w_{k}}).
\end{equation}

By Equations (1) and (2), we have $x_{V_{2}}=x_{U_{2}}.$ By Equations (3) and (4), we have $x_{V_{1}}=x_{U_{1}}.$ This completes the proof. \hfill$\square$

Let $u$ be a vertex of a connected $k$-uniform hypergraph $G$ with $|E(G)| \geq 1$. Let $G(u;p,q)$ be a $k$-unifom hypergraph obtained from $G$ by attaching two pendent paths $P=(u,e_1,u_1,\ldots,
\\u_{p-1},e_p,u_p)$ and $Q=(u,f_1,v_1,\ldots,v_{q-1},f_q,v_{q})$ at $u$.
Suppose that $G(u;p+1,q-1)$ is obtained from $G(u;p,q)$ by moving edge $f_q$ from $v_{q-1}$ to $u_p$. Then we say that $G(u;p+1,q-1)$ is obtained from $G(u;p,q)$ by grafting an edge.

The proof of the following Lemma 2.11 is similar to that of Xiao and Wang \cite{XW}.

\begin{lemma} \label{le:10}
Let $u$ be a vertex of a connected $k$-uniform hypergraph $G$ with $|E(G)| \geq 1$. If $p \geq q \geq 1$, then $q(G(u;p,q))>q(G(u;p+1,q-1))$.
\end{lemma}

\section{The supertree with the largest signless Laplacian spectral radius in $\mathbb{S}(m,d,k)$ }
\label{sec:ch-inco}

In this section,  we mainly determine the unique supertree with the largest signless Laplacian spectral radius among all supertrees in $\mathbb{S}(m,d,k)$ for $3\leq d\leq m-1$. And we determine the unique supertree with the second largest signless Laplacian spectral radius among all supertrees in $\mathbb{S}(m,d,k)$ for $d=3$.


Let $\mathbb{H}(m, d, k)$ be the set of $k$-uniform hypertrees with $m$ edges and diameter $d.$ Clearly, $\mathbb{H}(m, d, k) \subset \mathbb{S}(m, d, k).$ Let $S_{1}(m, d, k)$ be the $k$-uniform supertree obtained from a loose path $P:v_{1}, e_{1}, v_{2}, \ldots,v_{d}, e_{d}, v_{d+1}$ with length $d$ by attaching $m-d$ edges at vertex $v_{\lfloor\frac{d}{2}\rfloor+1}.$ By Definition 2.1, we know that $S_{1}(m, d, k)\in\mathbb{H}(m, d, k).$

Next we prove that $S_{1}(m, d, k)$ is the largest signless Laplacian spectral radius in $\mathbb{H}(m, d, k).$

\begin{theorem} \label{th:1}
Let $d\geq3$ and $m\geq d+1.$  Then $S_{1}(m, d, k)$ is the unique hypertree which attains the largest signless Laplacian spectral radius in $\mathbb{H}(m, d, k).$
\end{theorem}

\noindent\textbf{Proof.} Suppose $G$ is the hypertree which attains the largest signless Laplacian spectral radius in $\mathbb{H}(m, d, k)$. Let $x$ be the principal eigenvector of $\mathcal{Q}(G)$ corresponding to $q(G).$ Let $v_{1}, e_{1}, v_{2}, \ldots, e_{d}, v_{d+1}$ be the path with length $d$ of $G,$ where $e_{i}=\{v_{i}, v_{i,1}, ...,v_{i,k-2}, v_{i+1}\}$ for $1\leq i \leq d.$

\noindent\textbf{Claim 1.} If $d_{v_{i}}\geq 3$, $v_{i}$ only attaches pendent edges (except $e_{i-1}$ and $e_{i}$) for $2\leq i \leq d.$

Since the diameter of $G$ is $d$, we know that $v_{2}$ and $v_{d}$ attach pendent edges (except $e_{2}$ and $e_{d}$). Suppose that there exists a vertex $v_{i}$ that attaches a hypertree $H'$ of $G$, $3\leq i \leq d-1.$ Let the diameter of $H'$ be $d'$ and $u_{1}(=v_{i}),e'_{1},u_{2},...,e'_{d'},u_{d'+1}$ be the path with length $d'$ of $H'$. Without loss of generality, assume that $x_{u_{1}}\geq x_{u_{d'}}$. Then we obtain $G'$ from $G$ by moving an edge $e_{d'}$ from $u_{d'}$ to $u_{1}$. It is obvious that the diameter of $G'$ is $d$. By Lemma 2.4, we have $q(G')>q(G),$ a contradiction.

\noindent\textbf{Claim 2.} There exists a vertex $v_{i}$ attaching $m-d$ pendent edges, $2 \leq i \leq d$.

Suppose that $v_{j}, v_{j+1},\ldots, v_{j+s}$ attach $a_{j}, a_{j+1}, \ldots, a_{j+s}$ pendent edges of $G$, $2\leq j \leq j+s \leq d,$ respectively, where $a_{j}, a_{j+1}, \ldots, a_{j+s}\geq 0$. Assume that $x_{v_{i}}=\max\limits_{j \leq t \leq j+s}\{x_{v_{t}}\},$ then we obtain $G'$ from $G$ by moving $a_{t}$ pendent edges from $v_{t}$ to $v_{i}$, $j \leq t \leq j+s$ and $t\neq i.$ By Lemma 2.4, we have $q(G')>q(G),$ a contradiction. 

\noindent\textbf{Claim 3.} The vertex $v_{\lfloor\frac{d}{2}\rfloor+1}$ attaches $m-d$ pendent edges of $G$.

By Claims 1 and 2, we know that $G$ is the $k$-uniform hypertree obtained from a loose path $P:v_{1}, e_{1}, v_{2}, \ldots,v_{d}, e_{d}, v_{d+1}$ with length $d$ by attaching $m-d$ pendent edges at a vertex $v_{i}$, denoted $G(i),$~ $2 \leq i \leq d$. It is obvious that $G(i)\cong G(d+2-i).$ Here we prove that the spectral radius of $G(i)$ is greater than the spectral radius of $G(j)$, $2\leq j<i \leq \lfloor\frac{d}{2}\rfloor+1.$ 
Take two paths $P=v_{1},e_{1},v_{2},\ldots,v_{i}$ and $Q=v_{i},e_{i},v_{i+1},\ldots,v_{d+1}$ of $G(i)$. By Lemma 2.10, we know that the spectral radius of $G(i)$ is greater than the spectral radius of $G(j)$, $2\leq j<i \leq \lfloor\frac{d}{2}\rfloor+1.$ Thus we know that the vertex $v_{\lfloor\frac{d}{2}\rfloor+1}$ attaches $m-d$ pendent edges of $G$.




Hence, $G$ is obtained from a loose path $P:v_{1}, e_{1}, v_{2}, \ldots,v_{d}, e_{d}, v_{d+1}$ with length $d$ by attaching $m-d$ pendent edges at vertex $v_{\lfloor\frac{d}{2}\rfloor+1}.$ It is obviously that $G\cong S_{1}(m, d, k).$ Then $S_{1}(m, d, k)$ be the unique hypertree which attains the largest signless Laplacian spectral radius in $\mathbb{H}(m, d, k).$ This completes the proof. \hfill$\square$


\begin{theorem} \label{th:2}
Let $d\geq3$, $m\geq d+1$ and $G \in \mathbb{S}(m, d, k)\backslash \mathbb{H}(m, d, k).$ Then there exists a hypertree $G'\in\mathbb{H}(m, d, k)$ such that $q(G')>q(G).$
\end{theorem}

\noindent\textbf{Proof.} Let $v_{1}, e_{1}, v_{2}, \ldots, e_{d}, v_{d+1}$ be the path with length $d$ of $G,$ where $e_{i}=\{v_{i}, v_{i,1}, \ldots, v_{i,k-2},\\ v_{i+1}\}$ for $1 \leq i \leq d.$ Suppose that $e'_{1}, e'_{2}, \ldots, e'_{r}$ are all the branch edges in $E(G) \backslash \{e_{1}, e_{2}, \ldots, e_{d}\}.$ We obtain a hypergraph $G'$ from $G$ by edge-releasing on $e'_{j}$ at any vertex $v$ of $e'_{j}$ for $1 \leq j\leq r,$ and moving all the edges incident with one of vertices $v_{i,1}, \ldots, v_{i,k-2}$ from $v_{i,1}, \ldots, v_{i,k-2}$ to $v_{i}$ for $2 \leq i \leq d-1,$ respectively. Since there is no branch edge of $G',$ $G'$ is a hypertree. From Lemmas 2.6, 2.7, 2.8 and the diameter of $G'$ is $d,$ we have that $G'\in \mathbb{H}(m, d, k)$ and $q(G')>q(G).$ This completes the proof. \hfill$\square$

By Theorems 3.1, 3.2, we can obtain the following result.

\begin{theorem} \label{th:3}
Let $d\geq3$, $m\geq d+1$ and $G \in \mathbb{S}(m, d, k).$ Then $q(G)<q(S_{1}(m,d,k)).$
\end{theorem}

Let $S_{3}(m, d, k)$ be the $k$-uniform supertree obtained from a loose path $v_{1}, e_{1}, v_{2},\ldots,e_{d}, v_{d+1}$ by attaching $m-d$ edges at a pendent vertex $v_{0}$ of $e_{\lfloor\frac{d}{2}\rfloor+1}$ where $d$ is odd (see Figure 1).

\begin{center}
\includegraphics [width=7 cm, height=2 cm]{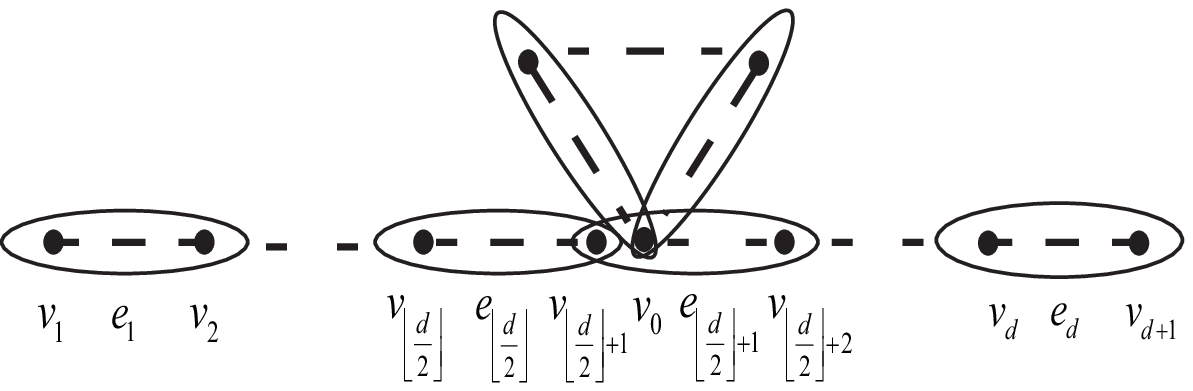}
\centerline{Figure 1: $S_{3}(m, d, k)$}
\end{center}

For $d=3$ and $m\geq d+2,$ let $S_{4}(m,d,k)$ be the $k$-uniform supertree obtained from a loose path $v_{1},e_{1},v_{2},e_{2},v_{3},e_{3},v_{4}$ by attaching $m-4$ pendent edges at vertex $v_{2}$ and attaching one pendent edge at vertex $v_{3}.$

\begin{theorem} \label{th:4}
(1) If $d=3$ and $m=d+1.$ Then we have $q(S_{1}(m, d, k))>q(S_{3}(m, d, k)).$
(2) If $d=3$ and $m\geq d+2.$ Then $q(G)<q(S_{4}(m, d, k))<q(S_{1}(m, d, k))$ for any $G \in \mathbb{S}(m, d, k)\backslash \{S_{1}(m, d, k), S_{4}(m, d, k)\}.$
\end{theorem}

\noindent\textbf{Proof.}
(1). For $d=3$ and $m=d+1,$ we know that $\mathbb{S}(m,d,k)=\{S_{1}(m,d,k),S_{3}(m,d,k)\}.$ By Theorem 3.3, we have $q(S_{1}(m, d, k))>q(S_{3}(m, d, k)).$

(2). Let $v_{1}, e_{1}, v_{2},e_{2}, v_{3},e_{3}, v_{4}$ be the path with length 3 of $G,$ where $e_{i}=\{v_{i}, v_{i,1}, \ldots, v_{i,k-2}\\, v_{i+1}\}$ for $1\leq i\leq 3.$ Since  $G \in \mathbb{S}(m, d, k),$ we know that $e_{1},e_{3}$ have $k-1$ pendent vertices. Let $x$ be the principal eigenvector of $\mathcal{Q}(G)$ corresponding to $q(G).$

{\bf Case 1.} Suppose that $e_{2}$ is a non-branch edge. Since $G\ncong S_{1}(m, d, k)$ and $G\ncong S_{4}(m, d, k)$, we know that $G$ respectively have $a$ pendent edges and $b$ pendent edges in $v_2$ and $v_3$ such that $a,~b\geq 3$. Without loss of generality, we assume that $x_{v_{2}}\geq x_{v_{3}}.$
We obtain a hypergraph $G'$ from $G$ by moving $b-2$ pendent edges from $v_{3}$ to $v_{2}$. It is obvious that $G'\cong S_{4}(m, d, k)$. By Lemma 2.4, we have $q(G)<q(G')= q(S_{4}(m, d, k)).$

{\bf Case 2.} Suppose that $e_{2}$ is a branch edge. Let $v_2, v_{2,i}, v_{2,i+1},\ldots, v_{2,i+s}, v_3$ attach $a, a_{i}, a_{i+1}, \\\ldots, a_{i+s}, b$ pendent edges of $G$ for $1\leq i \leq i+s \leq k-2,$ respectively, where $a_{i}, a_{i+1}, \ldots, a_{i+s}\geq0$ and $a,~b\geq 1$. 
Without loss of generality, we assume that $x_{v_{2,t}}=\max\limits_{i\leq j \leq i+s}\{x_{v_{2,j}}\},$ we obtain $G'$ from $G$ by moving $a_{j}$ edges from $v_{2,j}$ to $v_{2,t}$ of $G$ for all $i\leq j \neq t \leq i+s.$ Without loss of generality, we assume that $x_{v_{2}}=\max\{x_{v_{2}},x_{v_{2,t}},x_{v_{3}}\},$ we obtain $G''$ from $G'$ by moving $\sum_{j=i}^{i+s}a_{j}-1$ edges and $b-1$ edges from $v_{2,t}$ and $v_{3}$ to $v_{2}$ of $G$, respectively. It is obvious that $G''\cong S_{3}(m, d, k)$. By Lemma 2.4, we have $q(G)<q(G')<q(G'')=q(S_{3}(m, d, k)).$


Let $x'$ be the principal eigenvector of $\mathcal{Q}(S_{3}(m, d, k))$ corresponding to $q(S_{3}(m, d, k)).$ If $x'_{v_{2}}\geq x'_{v_{3}},$ we obtain $G'$ from $S_{3}(m, d, k)$ by moving $e_{3}$ from $v_{3}$ to $v_{2}.$ It is obvious that $G'\cong S_{4}(m, d, k)$. By Lemma 2.4, we have $q(S_{3}(m, d, k))<q(G')=q(S_{4}(m, d, k)).$
If $x'_{v_{2}}< x'_{v_{3}},$ we obtain $G''$ from $S_{3}(m, d, k)$ by moving $e_{1}$ from $v_{2}$ to $v_{3}.$ It is obvious that $G''\cong S_{4}(m, d, k)$. By Lemma 2.4, we have $q(S_{3}(m, d, k))<q(G'')=q(S_{4}(m, d, k)).$

Hence, if $d=3$ and $m\geq d+2,$ we have $q(G)<q(S_{4}(m, d, k))<q(S_{1}(m, d, k))$ for any $G \in \mathbb{S}(m, d, k)\backslash \{S_{1}(m, d, k), S_{4}(m, d, k)\}.$ This completes the proof. \hfill$\square$





Let $S_{2}(m, d, k)$ be the $k$-uniform supertree obtained from a loose path $v_{1}, e_{1}, v_{2}, \ldots,v_{d}, e_{d},\\ v_{d+1}$ by attaching $m-d$ edges at vertex $v_{\lfloor\frac{d}{2}\rfloor}$ (see Figure 2).

\begin{center}
\includegraphics [width=7 cm, height=2 cm]{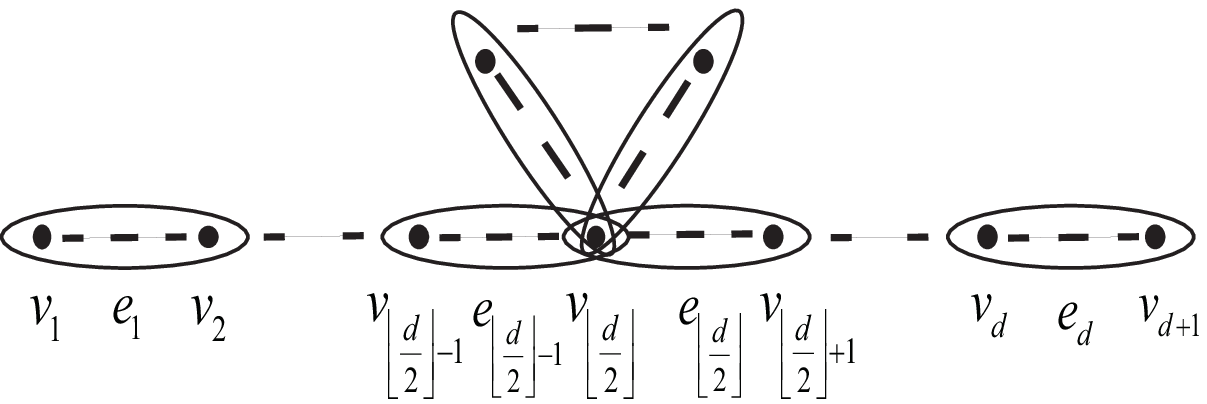}
\centerline{Figure 2: $S_{2}(m, d, k)$}
\end{center}

Xiao, Wang and Du \cite{XiWD} determined the first two largest spectral radii of uniform supertrees with given diameter, i.e., $\rho(G)<\rho(S_{2}(m,d,k))<\rho(S_{1}(m,d,k))$ for any $G\in\mathbb{S}(m,d,k)$. They compared $\rho(S_{2}(m, d, k))$ and $\rho(S_{3}(m, d, k))$ by using weighted incidence matrix. 
We may need to find the new methods and techniques to determine the second largest signless Laplacian spectral radius of $k$-uniform supertrees.

\begin{conjecture} \label{co:7}
Let $d\geq 4$ and $m\geq d+1.$ Then $S_{2}(m,d,k)$ has the second largest signless Laplacian spectral radius in $\mathbb{S}(m,d,k)$.
\end{conjecture}

\section{The supertree with the maximum signless Laplacian spectral radius with given number of pendant edges (vertices)}
In this section, by using a similar method of \cite{XW}, we determine the unique supertree with the maximum signless Laplacian spectral radius in $\mathbb{T}(n,p,k)$ and $\mathbb{G}(n,q,k)$, respectively.

Loose paths $P_{s_1},P_{s_2},£¬ \ldots, P_{s_q}$ are said to have almost equal lengths if $s_1, s_2, \ldots, s_q$ satisfy $|s_{i}-s_{j}| \leq 1$ for $1 \leq i \leq j \leq q$.

Let $T_{1}(n,p,k)$ be a $k$-uniform supertree on $n$ vertices obtained from a hyperstar $S_{p(k-1)+1,k}$ and $p$ loose paths of almost equal lengths by joining each pendant vertex of $S_{p(k-1)+1,k}$ to an end vertex of one loose path. In general, $T_{1}(n,\frac{n-1}{k-1},k) \cong S_{n,k}$, $T_{1}(n,1,k) \cong P_{k,k}$ and $T_{1}(n,2,k) \cong P_{n,k}$.

The proofs of the following results are similar to Lemma 4.1 and Theorem 4.2 of \cite{XW}, respectively. So we omit the proofs.

\begin{lemma} \label{le:4.1}
Let $e$ be an edge of a $k$-uniform hypergraph $G \in \mathbb{T}(n,p,k)$. If $e$ is a branch edge, then there exists a $k$-uniform supertree $G'\in \mathbb{T}(n,p,k)$ such that $q(G')>q(G)$.
\end{lemma}


\begin{theorem} \label{th:4.2}
 Let $G$ be a supertree in $\mathbb{T}(n,p,k)$, where $ 2 \leq p\leq \frac{n-1}{k-1}$. Then $q(G)<q(T_{1}(n,p,k))$, and the equality holds if and only if $G \cong T_{1}(n,p,k)$.
\end{theorem}





In \cite{XWL}, Xiao, Wang and Lu introduced the concept of BFS-supertree.

\begin{definition} \label{de:4.3}\cite{XWL}
Let $G=(V,E)$ be a $k$-uniform supertree with root $v_{0}$. A well-ordering $\prec$ of the vertices is called a breadth-first-search ordering (BFS-ordering for short) if all the following hold for all vertices:

$(a)$ $u\prec v$ implies $h(u)\leq h(v).$

$(b)$ $u\prec v$ implies $d_u\geq d_v.$

$(c)$ If ${u,u_{1}}\subset e_{1}\in E$ and ${v,v_{1}}\subset e_{2}\in E$ such that $u\prec v$, $h(u)=h(u_1)+1,$ $h(v)=h(v_1)+1,$ then $u_1\prec v_1.$

$(d)$ Suppose $u_1\prec u_2\prec\ldots\prec u_k$ for every edge $e=\{u_1,u_2,\ldots, u_k\}\in E,$ then there exists no vertex $v\in V\backslash e$ such that $u_i\prec v\prec u_{i+1},$ $2\leq i \leq k-1.$

\end{definition}

For a given degree sequence $\pi =(d_{0},d_{1},\ldots,d_{n-1}),$ let $\mathbb{S}_{\pi}$ be a set of $k$-uniform supertrees with  $\pi$  as its degree sequence.

\begin{lemma} \label{le:4.4}
Let $G=(V,E)$ be a $k$-uniform supertree attaining the maximum signless Laplacian spectral radius in $\mathbb{S}_{\pi}$ and $x$ be the principal eigenvector of $\mathcal{Q}(G).$ If $d_{u}>d_{v},$ then $x_{u}>x_{v}.$
\end{lemma}

\noindent\textbf{Proof.} The result can be proved by using a method similar to that used in Lemma 2.7 of \cite{XWL}. \hfill$\square$

\begin{lemma} \label{le:4.5}
Let $G=(V,E)$ be a $k$-uniform supertree attaining the maximum signless Laplacian spectral radius in $\mathbb{S}_{\pi}$ and $x$ be the principal eigenvector of $\mathcal{Q}(G).$ If $x_{u}\geq x_{v}$, then $d_{u}\geq d_{v}.$ Moreover, if $x_{u}= x_{v}$, then $d_{u}= d_{v}.$
\end{lemma}

\noindent\textbf{Proof.} The result can be proved by using a method similar to that used in Corollary 2.3 of \cite{XWL}. \hfill$\square$

\begin{theorem} \label{th:4.6}
For a given degree sequence $\pi=(d_{0},d_{1},\ldots,d_{n-1})$ of some $k$-uniform supertrees, if $G$ attains the maximum signless Laplacian spectral radius in $\mathbb{S}_{\pi}$, then $G$ is a BFS-supertree.
\end{theorem}

\noindent\textbf{Proof.} By Lemmas 2.4, 2.9, 4.4, 4.5 and Definition 4.3, the result can be proved by using a method similar to that used in Theorem 3.1 of \cite{XWL}. \hfill$\square$

\begin{proposition} \label{pro:4.7}\cite{XWL}
For a given degree sequence $\pi$ of some $k$-uniform supertrees, there exists a $k$-uniform supertree $G^{*}$ with degree sequence $\pi$ having BFS-ordering. Moreover, any two $k$-uniform supertrees with the same degree sequence and having BFS-ordering are isomorphic.
\end{proposition}

In \cite{XWL}, Xiao, Wang and Lu introduced a special BFS-supertree $G^{*}$.

\begin{theorem} \label{th:4.8}
For a given degree sequence $\pi$ of some $k$-uniform supertrees, BFS-supertree $G^{*}$ is a unique $k$-uniform supertree with maximum signless Laplacian spectral radius in $\mathbb{S}_{\pi}$.
\end{theorem}

\noindent\textbf{Proof.} Combining Definition 4.3, Theorem 4.6 and Proposition 4.7, we can obtain this result. \hfill$\square$

Let $G$ be a $k$-uniform supertree with $n$ vertices and $q$ pendent vertices for $n-\frac{n-1}{k-1}+1 \leq q \leq n$. If $q=n$, then $G \cong S_{k,k}$. And if $q=n-1$, then $G \cong S_{n,k}$.

The proof of the following result is similar to Lemma 6.1 in \cite{XW}.

\begin{lemma}\label{le:4.9}
 Let $G$ be a $k$-uniform supertree with the maximum signless Laplacian spectral radius in $\mathbb{G}(n,q,k)$ for $n-\frac{n-1}{k-1}+1 \leq q \leq n-2$. Then $G$ with degree sequence $\pi =(q+1+\frac{n-1}{k-1}-n,\underbrace{2,\ldots,2}\limits_{n-q-1},\underbrace{ 1,\ldots,1}\limits_{q})$.
\end{lemma}

By Theorem 4.8 and Lemma 4.9, we have the following theorem.

\begin{theorem}\label{th:4.10}
Let $G^{*}$ be the $k$-uniform BFS-supertree with degree sequence $\pi =(q+1+\frac{n-1}{k-1}-n,\underbrace{2,\ldots,2}\limits_{n-q-1},\underbrace{ 1,\ldots,1}\limits_{q})$, where $n-\frac{n-1}{k-1}+1 \leq q \leq n-1$. Then $G^{*}$ is the unique supertree attains the maximum signless Laplacian spectral radius in $\mathbb{G}(n,q,k)$.
\end{theorem}


\begin{thebibliography}{99}

\bibitem{BrHa} A.E. Brouwer, W.H. Haemers, Spectra of graphs, Universitext, Springer, New York, 2012. Available from: http://www.win.tue.nl/~aeb/2WF02/spectra.pdf.

\bibitem{CoDu} J. Cooper, A. Dutle, Spectra of uniform hypergraphs, Linear Algebra Appl. 436 (2012) 3268--3292.

\bibitem{CVRS} D. Cvertkovi$\acute{c}$, P. Rowlinson, S. Simi$\acute{c}$, An introduction to the theory of graph spectra, Cambridge University Press, Cambridge, 2010.

\bibitem{GuoS} J.M. Guo, J.Y. Shao, On the spectral radius of trees with fixed diameter, Linear Algebra Appl. 413 (2006) 131-147.

\bibitem{GuXC} S.G. Guo, G.H. Xu, Y.G. Chen, The spectral radius of trees with $n$ vertices and diameter $d$, Adv. Math. 6 (2005) 683-692.

\bibitem{HuQS} S.L. Hu, L.Q. Qi, J.Y. Shao, Cored hypergraphs, power hypergraphs and their Laplacian $H$-eigenvalues, Linear Algebra Appl. 439 (2013) 2980-2998.

\bibitem{KLS} L.Y. Kang, L.L. Liu, E.F. Shan, Sharp lower bounds for the spectral radius of uniform hypergraphs concerning degrees, The Electronic Journal of Combinatorics. 25(2) (2018) 1-13.

\bibitem{LiZB} H.F. Li, J. Zhou, C.J. Bu, Principal eigenvectors and spectral radii of uniform hypergraphs, Linear Algebra Appl. 544 (2018) 273-285.
\bibitem{LiSQ} H.H. Li, J.Y. Shao, L.Q. Qi, The extremal spectral radii of $k$-uniform supertrees, J. Comb. Optim. 32 (2016) 741-764.

\bibitem{LimL} L.H Lim, Singular values and eigenvalues of tensors: a variational approach. In; proceedings of the IEEE international workshop on computational advances in multi-sensor adaptive processing (CAMSAP 05). 1 (2005) 129-132.

\bibitem{Lim} L.H. Lim, Eigenvalues of tensors and some very basic spectral hypergraph theory, matrix computations and scientific computing seminar, April 16, 2008. http://www.atat.uchicago.edu/lekheng/work/mcsc2.

\bibitem{LMZW} H.Y. Lin, B. Mo, B. Zhou, W.M Weng, Sharp bounds for ordinary and signless Laplacian spectral radii of uniform hypergraphs, Appl. Math. Comput. 285 (2016) 217-227.

\bibitem{LiKY} L.L. Liu, L.Y. Kang, X.Y. Yuan, On the principal eigenvector of uniform hypergraphs, Linear Algebra Appl. 511 (2016) 430-446.

\bibitem{LuMa} L.Y. Lu, S.D. Man, Connected hypergraphs with small spectral radius, Linear Algebra Appl. 509 (2016) 206-227.

\bibitem{OYQY} C. Ouyang, L.Q. Qi, X.Y. Yuan, The first few unicyclic and bicyclic hypergraphs with largest spectral radii, Linear Algebra Appl. 527 (2017) 141-163.

\bibitem{Q} L.Q. Qi, Eigenvalues of a real supersymmetric tensor, J. Symbolic Comput. 40 (2005) 1302-1324.

\bibitem{Qi} L.Q. Qi, Symmetric nonnegative tensors and copositive tensors, Linear Algebra Appl. 439 (2013) 228-238.

\bibitem{QiL} L.Q. Qi, $H^{+}$-eigenvalues of Laplacian and signless Laplacian tensors, Commun. Math. Sci. 12 (2014) 1045-1064.

\bibitem{SKLS} L. Su, L.Y. Kang, H.H. Li, E.F. Shan, The matching polynomials and spectral radii of uniform supertrees. arXiv:1807.01180v1.

\bibitem{XW} P. Xiao, L.G. Wang, The maximum spectral radius of uniform hypergraphs with given number of pendant edges, Linear and Multilinear Algebra. 2018. DOI: 10.1080/03081087.2018.1453471.

\bibitem{XiWD} P. Xiao, L.G. Wang, Y.F. Du, The first two largest spectral radii of uniform supertrees with given diameter, Linear Algebra Appl. 536 (2018) 103-119.

\bibitem{XWL} P. Xiao, L.G. Wang, Y. Lu, The maximum spectral radii of uniform supertrees with given degree sequences, Linear Algebra Appl. 523 (2017) 33-45.

\bibitem{YuSS} X.Y. Yuan, J.Y. Shao, H.Y. Shan, Ordering of some uniform supertrees with largest spectral radii, Linear Algebra Appl. 495 (2016) 206-222.

\bibitem{YZLQ} J.J. Yue, L.P. Zhang, M. Lu, L.Q. Qi, The adjacency and signless Laplacian spectrai radius of cored hypergraphs and power hypergraphs, J. Oper. Res. Soc. China. 5 (2017) 27-43.

\end{thebibliography}
\end{document}